\documentclass[12pt,leqno]{article}
\usepackage{amsfonts,amsthm,amsmath}
\newtheorem{thm}{Theorem}

\theoremstyle{remark}

\newcommand{\FF}{\mathbb{F}}
\newcommand{\ZZ}{\mathbb{Z}}
\newcommand{\RR}{\mathbb{R}}

\newcommand{\ra}{\rangle}
\newcommand{\la}{\langle}

\begin{document}
\title{There exists no self-dual $[24,12,10]$ code
over $\FF_5$}

\author{
Masaaki Harada\thanks{
Department of Mathematical Sciences,
Yamagata University,
Yamagata 990--8560, Japan}
and 
Akihiro Munemasa\thanks{
Graduate School of Information Sciences,
Tohoku University,
Sendai 980--8579, Japan}
}

\maketitle

\begin{abstract}
Self-dual codes over $\FF_5$ exist for all even lengths.
The smallest length for which the largest minimum
weight among self-dual codes has not been determined
is $24$, and the largest minimum weight is either $9$ or $10$.
In this note, we show that there exists no self-dual 
$[24,12,10]$ code over $\FF_5$, 
using the classification of $24$-dimensional odd unimodular
lattices due to Borcherds.
\end{abstract}

\section{Introduction}

Let $\FF_5$ denote the finite field of order $5$.
An $[n,k]$ code $C$ over $\FF_5$ is a $k$-dimensional vector subspace
of $\FF_5^n$, where $n$ is called the length of $C$.
All codes in this note are codes over $\FF_5$.
An $[n,k,d]$ code is an $[n,k]$ code with minimum weight $d$.
A code $C$ is said to be {\em self-dual} if
$C=C^\perp$, where $C^\perp$ denotes
the dual code of $C$ under the standard inner product. 
A self-dual code of length $n$ exists if and only if 
$n$ is even.

As described in~\cite{RS-Handbook},
self-dual codes are an important class of linear codes for both
theoretical and practical reasons.
It is a fundamental problem to classify self-dual codes
of modest length 
and determine the largest minimum weight among self-dual codes
of that length.
Self-dual codes over $\FF_5$ were classified in \cite{LPS-GF5}
for lengths up to $12$.
The classification was extended to lengths $14$ and $16$
in \cite{HO-GF5}.
The largest minimum weights among self-dual codes of lengths $18,20$
and $22$ were determined in \cite{HO-GF5}, \cite{LPS-GF5} and 
\cite{HanKim}, respectively.
For length $24$, 
the largest minimum weight is either $9$ or $10$ \cite{LPS-GF5}.
In this note, we prove the following theorem.

\begin{thm}\label{thm}
There exists no self-dual $[24,12,10]$ code over $\FF_5$.
\end{thm}

Hence the largest minimum weight among
self-dual codes of length $24$ is exactly $9$.
The assertion of Theorem \ref{thm} was a question in
\cite[p.~192]{LPS-GF5}.

\section{Unimodular lattices and Construction A}

An $n$-dimensional (Euclidean) lattice $L$ is {\em unimodular} if
$L = L^{*}$, where
the dual lattice $L^{*}$ is defined as 
$L^{*} = \{ x \in {\RR}^n | \la x,y\ra \in \ZZ \text{ for all }
y \in L\}$ under the standard inner product $\la x, y\ra$.
The {\em norm} of a vector $x$ is $\la x, x\ra$.
The {\em minimum norm} of $L$ is the smallest 
norm among all nonzero vectors of $L$.
A unimodular lattice $L$ is {\em even} if all vectors of $L$ have even norms,
and {\em odd} if some vector has an odd norm.
The {\em kissing number} of $L$ is the number of vectors of minimum norm.

If $C$ is a self-dual code of length $n$, 
then 
\[
A_{5}(C) = \frac{1}{\sqrt{5}} 
\{x \in \ZZ^n \:|\: (x \bmod 5)\in C\}
\]
is an odd unimodular lattice,
where 
$(x \bmod 5)$ denotes
$(x_1 \bmod 5,\ldots,x_n \bmod 5)$ for $x=(x_1,x_2,\ldots,x_n)$.
This construction of lattices from codes
is called Construction A\@. 
If $C$ is a self-dual $[24,12,10]$ code over $\FF_5$, then
$A_5(C)$ is a $24$-dimensional odd unimodular 
lattice with minimum norm $\ge 2$.
The odd Leech lattice is a unique $24$-dimensional odd
unimodular lattice with minimum norm $3$.
There are $155$ non-isomorphic $24$-dimensional odd
unimodular lattices with minimum norm $2$ 
\cite{Bor} (see also \cite[Table 2.2]{SPLAG}).

\section{Proof}
In this section, we give a proof of Theorem \ref{thm}.

\begin{proof}
Let $C$ be a self-dual $[24,12,10]$ code over $\FF_5$.
As described in \cite{HanKim}, the Lee weight enumerator
(see \cite[p.~180]{LPS-GF5} for the definition) 
of $C$ is uniquely determined.
Since the coefficient of $x^{14}y^{10}$ in the
Lee weight enumerator is $528$,
$A_5(C)$ has minimum norm $2$ and kissing number $528$.
%
The only $24$-dimensional odd unimodular
lattice with minimum norm $2$ and kissing number $528$
is the $154$-th lattice
in~\cite[Table 17.1]{SPLAG},
which is the direct sum of two copies of the lattice $D_{12}^+$.
Thus
$A_5(C) = L_1 \oplus L_2$, where
for $i=1,2$, $L_i$ is isomorphic to $D_{12}^+$ when
restricted to the $12$-dimensional subspace $\mathbb{R}L_i$
of $\mathbb{R}^{24}$. In particular, both $L_1$ and $L_2$
have minimum norm $2$.

Let $e_i$ denote the unit vector
$(\delta_{i,1},\delta_{i,2},\ldots,\delta_{i,24})$
($i=1,2,\ldots,24$),
where $\delta_{i,j}$ is Kronecker's delta symbol.
We claim $\sqrt{5}e_i\in L_1$ or $\sqrt{5}e_i\in L_2$
for each $i\in \{1,\dots,24\}$.
Indeed, it suffices to
prove the claim for $i=1$. We may write
$\sqrt{5}e_1=a+b$, where $a\in L_1$ and $b\in L_2$.
Since the minimum norms of $L_1,L_2$ are both $2$,
$a\ne0$ and $b\ne0$ would imply
$\{\la a,a\ra,\la b, b\ra\}=\{2,3\}$.
We may assume without loss of generality that
$\la a,a\ra=2$, and write
$a=\frac{1}{\sqrt{5}}(c_1,\dots,c_{24})$. Then
$c_1=\la a,\sqrt{5}e_1 \ra =\la a,a+b \ra=2$
since $L_1=A_5(C)\cap L_2^\perp$, and hence
$10=5\la a,a \ra =\sum_{i=1}^{24}c_i^2=4+\sum_{i=2}^{24}c_i^2$.
This implies that the codeword $(c \bmod{5})\in C$
has weight less than $10$. This contradiction shows that
either $a=0$ or $b=0$, proving the claim.

Since the vectors $\sqrt{5}e_i$ ($i=1,2,\dots,24$) are
linearly independent and $\dim L_1=\dim L_2=12$, we may
assume without loss of generality that
\begin{align*}
&\sqrt{5}e_i\in L_1\quad\text{for }i=1,2,\dots,12,\text{ and}\\
&\sqrt{5}e_i\in L_2\quad\text{for }i=13,14,\dots,24.
\end{align*}
Then
\begin{align*}
L_1&=A_5(C)\cap L_2^\perp
=A_5(C)\cap \bigoplus_{i=1}^{12}\frac{1}{\sqrt{5}}\ZZ e_i,\\
L_2&=A_5(C)\cap L_1^\perp
=A_5(C)\cap \bigoplus_{i=13}^{24}\frac{1}{\sqrt{5}}\ZZ e_i.
\end{align*}
Define codes $C_1,C_2$ by
\begin{align*}
C_1&=\{((c_1,\dots,c_{12})\bmod{5})\mid
\frac{1}{\sqrt{5}}(c_1,\dots,c_{12},0,\dots,0)\in  
L_1\},\\
C_2&=\{((c_{13},\dots,c_{24})\bmod{5})\mid
\frac{1}{\sqrt{5}}(0,\dots,0,c_{13},\dots,c_{24})\in  
L_2\}.
\end{align*}
Then for $c=(c_1,c_2,\dots,c_{24})\in\ZZ^{24}$, we have
\begin{align*}
&(c\bmod{5})\in C
\\ &\iff
\frac{1}{\sqrt{5}}c\in L_1\oplus L_2
\\ &\iff
c=a_1+a_2\text{ for some }a_1,a_2\in\ZZ^{24}\text{ with }
\frac{1}{\sqrt{5}}a_1\in L_1,\;
\frac{1}{\sqrt{5}}a_2\in L_2
\\ &\iff
\frac{1}{\sqrt{5}}(c_1,\dots,c_{12},0,\dots,0)\in L_1\text{ and }
\frac{1}{\sqrt{5}}(0,\dots,0,c_{13},\dots,c_{24})\in L_2
\\ &\iff
((c_1,\dots,c_{12})\bmod{5})\in C_1\text{ and }
((c_{13},\dots,c_{24})\bmod{5})\in C_2.
\end{align*}
Hence $C$ is decomposable into the direct sum of the two codes
$C_1,C_2$, each of which is of length $12$.
Since $(C_1 \oplus C_2)^\perp = C_1^\perp \oplus C_2^\perp$
(see \cite[Exercise 30]{Huffman-Pless}) and $C$ is self-dual, both  
$C_1$ and $C_2$ are self-dual.
However,
no self-dual code of length $12$ has
minimum weight $\ge 10$. This is a contradiction, and
the proof is complete.
\end{proof}



\begin{thebibliography}{99}
\bibitem{Bor} R.E. Borcherds,
The Leech lattice and other lattices,
Ph.D. Dissertation, Univ.\ of Cambridge, 1984.

\bibitem{SPLAG} J.H. Conway and N.J.A. Sloane,
{\sl Sphere Packing, Lattices and Groups (3rd ed.)},
Springer-Verlag, New York, 1999.

\bibitem{HanKim} S. Han and J.-L. Kim,
On self-dual codes over $\Bbb F\sb 5$,
{\sl Des.\ Codes Cryptogr.}
{\bf 48} (2008), 43--58. 


\bibitem{HO-GF5} M. Harada and P.R.J. \"Osterg\aa rd,
On the classification of self-dual codes over~$\FF_5$,
{\sl Graphs Combin.}
{\bf 19} (2003), 203--214. 

\bibitem{Huffman-Pless} W.C. Huffman and V. Pless, 
{\sl Fundamentals of Error-Correcting Codes}, 
Cambridge University Press, Cambridge, 2003. 

\bibitem{LPS-GF5} J.S. Leon, V. Pless and N.J.A. Sloane,
{Self-dual codes over {\rm GF$(5)$}},
{\sl J.\ Combin.\ Theory~Ser.~A\/}
{\bf 32} (1982) 178--194.

\bibitem{RS-Handbook} E.\ Rains and N.J.A.\ Sloane,
{``Self-dual codes,''} {Handbook of Coding Theory},
V.S. Pless and W.C. Huffman (Editors),
Elsevier, Amsterdam 1998, pp.\ 177--294.

\end{thebibliography}
\end{document}